\documentclass[11pt,fleqn,twoside]{article}
\usepackage{amsfonts,amssymb,latexsym}
\makeatletter
\newcommand{\prava}[1]{\small\it
\begin{flushleft}
Copyright \copyright \ 1999 by  #1
\end{flushleft}}

\newcommand{\name}[1]{\begin{flushleft}
                       \LARGE \bf #1
                       \end{flushleft}\vspace{-3mm}}

\newcommand{\Author}[1]{\begin{flushleft}
                       \it #1 \end{flushleft}}

\newcommand{\Adress}[1]{\begin{flushleft}
                       \it #1 \end{flushleft}}

\newcommand{\Date}[1]{\begin{flushleft}
                      \small  \it #1 \end{flushleft}}

\newcommand{\ehkol}{Author \ name}
\newcommand{\ohkol}{Article \ name}
\renewcommand{\@evenhead}{
\hspace*{-3pt}\raisebox{-15pt}[\headheight][0pt]{\vbox{\hbox to \textwidth 
{\thepage \hfil \ehkol}\vskip4pt \hrule}}}
\renewcommand{\@oddhead}{
\hspace*{-3pt}\raisebox{-15pt}[\headheight][0pt]{\vbox{\hbox to \textwidth 
{\ohkol \hfil \thepage}\vskip4pt\hrule}}}
\renewcommand{\@evenfoot}{}
\renewcommand{\@oddfoot}{}

\newcommand{\be}{\begin{equation}}
\newcommand{\ee}{\end{equation}}
\newcommand{\ba}{\hspace*{-5pt}\begin{array}}
\newcommand{\ea}{\end{array}}

\newcommand{\ds}{\displaystyle}
\makeatother

\begin{document}

\thispagestyle{empty}
\setcounter{page}{246}
\renewcommand{\ehkol}{A.K. Prykarpatsky and D. Blackmore}
\renewcommand{\ohkol}{Versal Deformations of a Dirac Type
Dif\/ferential Operator}

\begin{flushleft}
\footnotesize \sf
Journal of Nonlinear Mathematical Physics \qquad 1999, V.6, N~3,
\pageref{pryk-black_fp}--\pageref{pryk-black_lp}.
\hfill {\sc Letter}
\end{flushleft}

\vspace{-5mm}

\renewcommand{\footnoterule}{}
{\renewcommand{\thefootnote}{}
 \footnote{\prava{A.K. Prykarpatsky and D. Blackmore}}}

\name{Versal Deformations of a Dirac Type \\
Dif\/ferential Operator}\label{pryk-black_fp}

\Author{Anatoliy K. PRYKARPATSKY~$^{\dag}$ and Denis BLACKMORE~$^\ddag$}

\Adress{$\dag$~Department of Applied
Mathematics at AGH, Cracow 30-059, Poland;\\[1.8mm]
~~Department of Nonlinear Mathematical Analysis at IAPMM of NAS,\\
~~Lviv 290601, Ukraina\\
~~E-mail: prika@mat.agh.edu.pl\\[2mm]
$\ddag$~Department of Mathematical Sciences and Center for Applied
Mathematics \\
~~and Statistics, New Jersey Institute of Technology, Newark, NJ
07102-1982, USA\\
~~E-mail: deblac@chaos.njit.edu}

\Date{Received November 10, 1998; Revised February 25, 1999; Accepted
April 1, 1999}

\begin{abstract}
\noindent
If we are given a smooth dif\/ferential operator in the variable $x\in
{\mathbb R}/2\pi {\mathbb Z},$ its normal form, as is well known, is
the simplest form obtainable by means of the $\mbox{Dif\/f}(S^1)$-group action
on the space of all such operators. A versal deformation of this
operator is a normal form for some parametric inf\/initesimal family
including the operator. Our study is devoted to analysis of versal
deformations of a Dirac type dif\/ferential operator using the theory
of induced $\mbox{Dif\/f}(S^1)$-actions endowed with centrally
extended Lie-Poisson brackets. After constructing a general
expression for tranversal deformations of a
Dirac type dif\/ferential operator, we interpret it via the Lie-algebraic
theory of induced $\mbox{Dif\/f}(S^1)$-actions on a special Poisson manifold and
determine its generic moment mapping. Using a Marsden-Weinstein reduction
with respect to certain Casimir generated distributions, we describe a wide
class of versally deformed Dirac type dif\/ferential operators depending on
complex parameters.
\end{abstract}

\renewcommand{\theequation}{\arabic{section}.\arabic{equation}}

\section{Introduction}

Suppose we are given the linear 2-vector f\/irst order Dirac dif\/ferential
operator on the real axis ${\mathbb R}$:
\begin{equation}
L_{{}\lambda }f:=-\frac{df}{dx}+l_{{}\lambda }[u,v;z]f,\qquad l_{{}\lambda
}[u,v;z]:=\left(
\begin{array}{cc}
z-\lambda &  u \\
 v & \lambda -z
\end{array}
\right)
\end{equation}
acting on the Sobolev space $W_{2,loc}^{(1)}({\mathbb R};{\mathbb C}^2)$
and depending on $2\pi$-periodic coef\/f\/icients $u,v,z\in C^\infty
({\mathbb R}/2\pi {\mathbb Z};{\mathbb C})$ and a complex parameter
$\lambda \in {\mathbb C}.$ The variety of all operators (1.1),
parametrized by $\lambda ,$ will be denoted by $\mathcal{L}_{{}\lambda }.$

Let $\mathcal{A}:=\mbox{Dif\/f}(S^1)$ be the group of orientation
preserving dif\/feomorphisms of the circle~$S^1.$ A group action
of $\mathcal{A}$ on $\mathcal{L}_{{}\lambda }$ can be def\/ined as
follows: Fixing a parametrization of $S^1,$ i.e., a $C^\infty $
covering $p:{\mathbb R} \rightarrow S^1$ such that the mapping
$p:[a,a+2\pi )\rightleftharpoons S^1$ is one-to-one for every real
$a$ and $p(x+2\pi )=p(x)$ for all $x\in {\mathbb R},$ each $\phi
\in \mathcal{A}$ can obviously be represented by a smooth mapping
$\phi :{\mathbb R}\rightarrow {\mathbb R}$ such that
\begin{equation}
\phi (\xi +2\pi )=\phi (\xi )+2\pi \qquad \mbox{and} \qquad
\phi ^{\prime }(\xi )>0
\end{equation}
for all $\xi \in {\mathbb R}.$ Upon making the change of variables
\begin{equation}
x=\phi (\xi ), \qquad f\left( \phi (\xi )\right) =\Phi (\xi
)\tilde{f}(\xi ),
\end{equation}
with $\phi \in \mathcal{A},$ $\Phi \in G:=C^\infty \left( {\mathbb
R}/2\pi {\mathbb Z};SL(2;{\mathbb C})\right) $ and $x,\xi \in
{\mathbb R},$ in~(1.1), it is easy to
see that the dif\/ferential operator $L_{{}\lambda }$ transforms into
$L_{{}\lambda }^{(\phi ,\Phi )}:W_2^{(1)}\rightarrow W_2^{(1)}$
def\/ined as
\begin{equation}
L_{{}\lambda }^{(\phi ,\Phi )}\tilde{f}(\xi ):=-\frac{d\tilde{f}}{d\xi }
+l_{{}\lambda }^{(\phi ,\Phi )}[u,v;z]\tilde{f},
\end{equation}
where
\begin{equation}
l_{{}\lambda }^{(\phi ,\Phi )}[u,v;z]:=-\Phi ^{-1}(\xi )\frac{d\Phi
(\xi )}{d\xi }+\phi ^{\prime }(\xi )\Phi ^{-1}(\xi )l_{{}\lambda
}[u,v;z]\Phi (\xi ).
\end{equation}

We assume now that the matrix $\Phi (\xi )$ is chosen so that
$l_{{}\lambda}^{(\phi ,\Phi )}[u,v;z]=l_{{}\lambda
}[\tilde{u},\tilde{v};\tilde{z}]$ for all $\lambda \in {\mathbb C}$
and some mapping $(\tilde{u},\tilde{v};\tilde{z})^T\in C^\infty
({\mathbb R}/2\pi {\mathbb Z};{\mathbb C}^2\times {\mathbb C}).$
Whence we obtain an induced nonlinear transformation $A^{*}(\phi
,\Phi ):\mathcal{L}_{{}\lambda }\rightarrow \mathcal{L}_{{}\lambda
},$ $(\phi ,\Phi )\in \mathcal{A\times }G,$ where
\begin{equation}
A^{*}(\phi ,\Phi )l_{{}\lambda }[u,v;z]:=
l_{{}\lambda}^{(\phi,\Phi)}[u,v;z]
\end{equation}
for all mappings in $C^\infty ({\mathbb R}/2\pi {\mathbb Z};{\mathbb
C}^2\times {\mathbb C}).$ This together with expression~(1.5)
determines an automorphism $A^{*}$ of $\mathcal{A},$ for a f\/ixed
$\Phi ,$ that we shall study in detail. We are primarily interested
in describing normal forms and versal deformations of~(1.1) with
respect to the automorphism $A^{*}.$

As is well known (see~[1, 2, 5]), a normal form of the operator~(1.1)
is the simplest (in some sense) representative of its orbit under the
group action of $\mathcal{A}$ on the space $\mathcal{L}_{{}\lambda }.$ A
versal deformation of~(1.1) is a normal form for a stable parametric
inf\/initesimal family including~(1.1). As will be shown below, all such
deformations can be described by means of Lie-algebraic analysis of this
group action on $\mathcal{L}_{{}\lambda }$ and an associated momentum
mapping reduced on certain invariant subspaces.

\setcounter{equation}{0}

\section{Lie-algebraic structure of the $\mathcal{A}$-action}

Let us consider the loop group $G:=G_{S^1}\left( SL(2;{\mathbb
C})\right) $ of all smooth mappings $S^1\rightarrow SL(2;{\mathbb
C})$ and its corresponding group $\mathcal{A}$-action on a functional
manifold $M\subset C^\infty ({\mathbb R}/2\pi {\mathbb Z};{\mathbb
C}^3),$ which is assumed to be equivariant; that is, the diagram
\begin{equation}
\ba{rcccl} & M & \stackrel{l}{\rightarrow } & \mathcal{G}^{*}  &\\
A_{{}\Phi }\!\!\! \!\! &\downarrow &  & \downarrow & \!\!\!\!\!
Ad_{{}\Phi ^{-1}}^{*}
\\ & M & \stackrel{l}{\rightarrow } & \mathcal{G}^{*} & \ea
\end{equation}
commutes for all $l$ in the adjoint $\mathcal{G}^{*}$ of the loop Lie
algebra and $\Phi \in G.$ Whence we can def\/ine on $M$ a natural Poisson
structure that induces the following canonical Lie-Poisson structure
on $\mathcal{G}^{*}$: for any $\gamma ,\mu \in D(\mathcal{G}^{*}),$
\begin{equation}
\{\gamma ,\mu \}:=(l,[\nabla \gamma (l),\nabla \mu (l)]).
\end{equation}
Here $(\cdot ,\cdot )$ is the usual Killing type nondegenerate, symmetric,
invariant scalar product on the loop Lie algebra
$\mathcal{G}=C_{S^1}(sl(2;{\mathbb C}))$, i.e. for any $a,b\in
\mathcal{G},$
\begin{equation}
(a,b):=\int_0^{2\pi }dx\; Sp(ab)
\end{equation}
and $\nabla :D(\mathcal{G}^{*})\rightarrow \mathcal{G}$ is def\/ined as $
(\nabla \gamma (l),\delta l):=\frac d{d\epsilon }\gamma (l+\epsilon \delta
l)\mid _{{}\epsilon =0}$ for any $\delta l\in \mathcal{G}^{*},\gamma
\in D(\mathcal{G}^{*}).$

In order to address the problems posed in Section 1, we need to centrally
extend the group action $A_{{}\Phi }:M\rightarrow M,$ $\Phi \in G,$ as
follows: for $\hat{\Phi}:=(\Phi ,c)\in \hat{G}:=G\times {\mathbb C}$ the
corresponding action $A_{{}\hat{\Phi}}:M\rightarrow M$ is def\/ined so that
the diagram
\begin{equation}
\ba{rcccl} & M & \stackrel{\hat{l}}{\rightarrow } &
\hat\mathcal{G}{}^{*} &\\ A_{{}\hat{\Phi}}\!\!\!\!\! & \downarrow
& & \downarrow & \!\!\!\!\! Ad_{{}\hat{\Phi}^{-1}}^{*} \\ & M &
\stackrel{\hat{l}}{\rightarrow } & \hat\mathcal{G}{}^{*}& \ea
\end{equation}
commutes for all $\hat{\Phi}\in \hat{G}$ and $\hat{l}=(l,c)\in
\hat\mathcal{G}{}^{*}.$ This leads to the following
(unique!) choice of the extended $Ad^{*}$-action in~(2.4):
\begin{equation}
Ad_{{}\hat{\Phi}^{-1}}^{*}:(l,c)\in \mathcal{G}^{*}\rightarrow \left( \phi
^{\prime }(\xi )Ad_{{}\Phi ^{-1}}l(x)-c\Phi ^{-1}\frac{d\Phi }{d\xi },c\right)
\end{equation}
for all $\hat{\Phi}\in \hat{G},$ $l\in \mathcal{G}^{*}$ at $\xi \in {\mathbb R},$
$x=\phi (\xi )$ and $c\in {\mathbb C}.$ This expression follows from the fact
that the loop Lie algebra $\mathcal{G}$ admits only the central extension
$\hat\mathcal{G}\oplus {\mathbb C}.$ As the homology groups
$H^1(\mathcal{G})=0$ and $H^2(\mathcal{G})=1,$ it is represented as
\begin{equation}
\left[ (a,\alpha ),(b,\beta )\right] :=\left( [(a,b)],(a,db/dx)\right)
\end{equation}
for any $a,b\in \mathcal{G}$ and $\alpha ,\beta \in {\mathbb C}.$ Taking $c$ to
be unity and def\/ining an appropriate di\-f\/feo\-mor\-phism $x\rightarrow \phi
(x)=\xi $ of ${\mathbb R},$ it is easy to see that $Ad_{{}\hat{\Phi}^{-1}}^{*}$
has the same structure element as that of the action $A^{*}(\phi ,\Phi )$ on
$\mathcal{L}_{\lambda}$ def\/ined above. Whence it is clear that our
Lie-algebraic analysis is intimately connected with the structure of the
$G$-orbits induced by the dif\/feomorphism group $\mathcal{A}=\mbox{Dif\/f}(S^1).$

We def\/ine a natural Lie-Poisson bracket on the adjoint space
$\hat\mathcal{G}{}^{*}$ as follows: for any $\gamma ,\mu
\in D(\hat\mathcal{G})\subset \hat\mathcal{G}{}^{*},$
\begin{equation}
\{\gamma ,\mu \}_0:=\left( l,[\nabla \gamma (l),\nabla \mu (l)]\right)
+\left( \nabla \gamma (l),\frac{d\nabla \mu (l)}{dx}\right) ,
\end{equation}
and deform it into a brackets pencil using a constant parameter $\lambda \in
{\mathbb C}$ via
\begin{equation}
\{\gamma ,\mu \}_0\stackrel{\lambda }{\rightarrow }\{\gamma ,\mu
\}_{{}\lambda }:=(\nabla \gamma (l),\frac d{dx}\nabla \mu (l))+(l+\lambda
J,[\nabla \gamma (l),\nabla \mu (l)]),
\end{equation}
where $J\in sl^{*}(2;{\mathbb C})$ is chosen here to be the constant matrix
\begin{equation}
J=\left(
\begin{array}{cc}
1 & 0 \\
0 & -1
\end{array}
\right) .
\end{equation}
The following compatibility condition is almost obvious [8, 10].

\medskip

\noindent
{\bf Lemma 2.1.} {\it A pencil of brackets (2.8) is a Poisson
brackets pencil for each $\lambda \in {\mathbb C}$ and $J\in sl^{*}(2;{\mathbb C}),$
i.e. it is compatible.}

\medskip

\noindent
{\bf Proof.} It is well known that the Lie derivative of a
Poisson bracket is also a Poisson bracket if and only if
\begin{equation}
\{\gamma ,\mu \}_1:=\mathfrak{L}_K\{\gamma ,\mu \}_0-\{\mathfrak{L}_K\gamma ,\mu
\}_0-\{\gamma ,\mathfrak{L}_K\mu \}_0
\end{equation}
satisf\/ies the Jacobi identity for all $\gamma ,\mu \in D(\mathcal{G}^{*}),$
where $\mathfrak{L}_K$ is the Lie derivative with respect to a vector f\/ield $K:
\mathcal{G}^{*}\rightarrow T(\mathcal{G}^{*}).$ Choosing $K(l):=J,$ it is
easy to verify that the bracket (2.10) satisf\/ies the Jacobi identity and is
the usual Poisson bracket on $\mathcal{G}^{*}.$ Consequently, the Poisson
bracket (2.10) is also a Poisson bracket along a generic orbit of the vector
f\/ield $dl/d\lambda =J,$ hence the deformation (2.8) is also Poisson, as was
to be proved.

\setcounter{equation}{0}

\section{Casimir functionals and reduction problem}

A Casimir functional $h\in I_{{}\lambda }(\hat\mathcal{G}{}^{*})$ is
def\/ined, as usual, as a functional $h\in D(\hat\mathcal{G}{}^{*})$ that is
invariant with respect to the following $\lambda$-deformed $Ad_{{}\hat{\Phi}^{-1}}^{*}$-action:
\begin{equation}
Ad_{{}\hat{\Phi}^{-1}}^{*}:(l,1)\in \hat\mathcal{G}{}^{*}\rightarrow
\left( Ad_{{}\Phi ^{-1}}^{*}(l+\lambda J)-\Phi ^{-1}\frac{d\Phi}{dx},1\right)
\end{equation}
for any $\Phi \in G,$ $l\in \mathcal{G}^{*}$ and $\lambda \in {\mathbb C}.$ It
is easy to see from this def\/inition that
$h\in I_{{}\lambda }(\hat\mathcal{G}{}^{*})$ if the equation
\begin{equation}
\frac{d\nabla h(l)}{dx}=\left[ l+\lambda J,\nabla h(l)\right]
\end{equation}
is satisf\/ied for all $\lambda \in {\mathbb C}.$ Assuming further that there
exists an asymptotic expansion of the form
\begin{equation}
h(\lambda )\sim \sum\limits_{j\in {\mathbb Z}_{+}}h_j\lambda ^{-j}
\end{equation}
as $\left| \lambda \right| \rightarrow \infty$, one can readily verify that
$h_0\in I_1(\hat\mathcal{G}{}^{*})$ and that for all
$j,k\in {\mathbb Z}_{+}$
\be
\{h_j,h_k\}_0 =0=\{h_j,h_k\}_1,  \qquad
 \{\gamma ,h_j\}_0
=\{\gamma ,h_{j+1}\}_1,
\ee
where $\gamma \in
D(\hat\mathcal{G}{}^{*})$ is arbitrary.

Let us now consider the action (2.1) at a f\/ixed
$l=l[u,v;z]\in \hat\mathcal{G}^{*}.$ It is easy to see that this action does not
necessarily preserve the form of the element $l.$ Thus we must reduce the
initial $\hat{G}$-action on $\hat\mathcal{G}^{*}$ to an
appropriate subgroup; for this we develop the reduction procedure employed
in [8--10].

Def\/ine the distribution
\begin{equation}
D_1:=\left\{ K\in T(\hat\mathcal{G}^{*}):K(l)=[J,\nabla
\gamma (l)],l\in \hat\mathcal{G}^{*},\gamma \in
D(\hat\mathcal{G}^{*})\right\} .
\end{equation}
$D_1$ is integrable, that is $[D_1,D_1]\subset D_1,$ since the bracket $\{\cdot ,\cdot \}_1$
is Poisson. Now def\/ine another distribution
\begin{equation}
D_0:=\left\{ K\in T(\hat\mathcal{G}^{*}):K(l)=[l-\frac
d{dx},\nabla h_0],h_0\in I_1(\hat\mathcal{G}^{*})\right\} ,
\end{equation}
which is clearly also integrable on $\hat\mathcal{G}^{*},$ since
$[D_0,D_0]\subset D_0.$ The set of maximal integral submanifolds of~(3.6)
generates the foliation $\hat\mathcal{G}_J^{*}\backslash D_0$
whose leaves are the intersections of f\/ixed integral
submanifolds
$\hat\mathcal{G}_J^{*}\subset \hat\mathcal{G}^{*}$
of the distribution $D_1$ passing through an
element $l[u,v;z]\in \hat\mathcal{G}^{*}.$ If the
foliation $\hat\mathcal{G}_J^{*}\backslash D_0$ is
suf\/f\/iciently smooth, one can def\/ine the quotient manifold
$\hat\mathcal{G}_{\mbox{\scriptsize red}}^{*}:=
\hat\mathcal{G}_J^{*}/(\hat\mathcal{G}_J^{*}\backslash
D_0)$ with its associated projection mapping $\hat\mathcal{G}_J^{*}
\rightarrow \hat\mathcal{G}_{\mbox{\scriptsize red}}^{*}.$ To continue
this line of reasoning, we shall obtain explicit constructions of the
objects introduced.

$D_1$ is obviously generated by the vector f\/ields
\begin{equation}
\frac{dl}{dt}=\left(
\begin{array}{cc}
0 & 2b \\
-2c & 0
\end{array}
\right) ,\qquad \nabla \gamma (l)=\left(
\begin{array}{cc}
a & b \\
c & -a
\end{array}
\right) ,
\end{equation}
where $t$ is a complex evolution parameter and $l\in \hat\mathcal{G}_J^{*},$ where
$\hat\mathcal{G}_J^{*}\subset \hat\mathcal{G}^{*}$
is the isotropy Lie subalgebra of the element $J\in \hat\mathcal{G}^{*}.$
Hence the integral submanifold $\hat\mathcal{G}_J^{*}$
consists of orbits of an element $l=l[u,v;z]\in \hat\mathcal{G}^{*},$ with
$z\in {\mathbb C},$ with respect to the vector f\/ields~(3.7). The distribution $D_0$ on
$T(\hat\mathcal{G}^{*})$ is generated by the vector f\/ields
\begin{equation}
\frac{dl}{d\tau }=\left(
\begin{array}{cc}
-\chi _x & -2u\chi \\
2v\chi & \chi _x
\end{array}
\right) ,\qquad \nabla h_0(l)=\left(
\begin{array}{cc}
\chi & 0 \\
0 & -\chi
\end{array}
\right) ,
\end{equation}
where $\tau $ is a complex evolution parameter and
$l=l[u,v;z]\in \hat\mathcal{G}^{*}.$

It follows immediately from (3.8) that
\begin{equation}
\frac{dz}{d\tau }=-\chi _x,\qquad \frac{du}{d\tau }=-2u\chi \qquad \mbox{and}\qquad
\frac{dv}{d\tau }=2v\chi
\end{equation}
for all $\tau \in {\mathbb R}$ along $D_0.$ Eliminating the variable $\chi $
from (3.9), we obtain
\begin{equation}
\frac d{d\tau }\left[ \frac d{dx}(\ln u)-2z\right] =0=\frac d{d\tau }\left[
\frac d{dx}(\ln v)+2z\right] ;
\end{equation}
that is, the mapping
\begin{equation}
\hat\mathcal{G}^{*}\ni l=\left(
\begin{array}{cc}
z & u \\
v & -z
\end{array}
\right) :\stackrel{\nu }{\rightarrow }\left(
\begin{array}{cc}
 0 & \exp (\partial ^{-1}\alpha ) \\
\exp (\partial ^{-1}\beta ) & 0
\end{array}
\right) \rightarrow \hat\mathcal{G}_{\mbox{\scriptsize red}}^{*},
\end{equation}
where
\begin{equation}
\alpha :=u_xu^{-1}-2z,\qquad \beta :=v_xv^{-1}+2z,
\end{equation}
explicitly determines the reduction $\nu :\hat\mathcal{G}^{*}\rightarrow
\hat\mathcal{G}_{\mbox{\scriptsize red}}^{*}$ discussed
above. We are now in a position to compute the bracket~(2.8) reduced upon
the submanifold $\hat\mathcal{G}_{\mbox{\scriptsize red}}^{*}$ by
def\/ining the functionals $\lambda ,\mu \in D(\hat\mathcal{G}^{*})$ to be
constant along the distribution $D_0,$ that is
\begin{equation}
\gamma :=\tilde{\gamma}\circ \nu ,\qquad \mu :=\tilde{\mu}\circ \nu ,
\end{equation}
for any $\tilde{\gamma},\tilde{\mu}\in D(\hat\mathcal{G}_{\mbox{\scriptsize red}}^{*}).$
From~(3.12) one readily obtains the expressions
\be
\ba{l}
\left. \nabla \gamma (l)\right| _{l\in \hat\mathcal{G}_{\mbox{\scriptsize red}}^{*}} =
\left(
\begin{array}{cc}
\ds \frac{\delta \tilde{\gamma}}{\delta \beta }-\frac{\delta \tilde{\gamma}}{
\delta \alpha } & \ds -\frac 1v\left( \frac{\delta \tilde{\gamma}}{\delta \beta }
\right) _x
\vspace{3mm}\\
\ds -\frac 1u\left( \frac{\delta \tilde{\gamma}}{\delta \alpha }\right) _x &
\ds \frac{\delta \tilde{\gamma}}{\delta \alpha }-\frac{\delta \tilde{\gamma}}{
\delta \beta }
\end{array}
\right) ,
\vspace{3mm}\\
\left. \nabla \mu (l)\right| _{l\in \hat\mathcal{G}_{\mbox{\scriptsize red}}^{*}} =\left(
\begin{array}{cc}
\ds \frac{\delta \tilde{\mu}}{\delta \beta }-\frac{\delta \tilde{\mu}}{\delta \alpha } &
\ds -\frac 1v\left( \frac{\delta \tilde{\mu}}{\delta \beta }\right) _x
\vspace{3mm}\\
\ds -\frac 1u\left( \frac{\delta \tilde{\mu}}{\delta \alpha }\right) _x &
\ds \frac{\delta \tilde{\mu}}{\delta \alpha }-\frac{\delta \tilde{\mu}}{\delta \beta }
\end{array}
\right) ,
\ea
\ee
which satisfy the desired identities
\begin{equation}
\left( \nabla \gamma (l),dl/d\tau \right) =0=\left( \nabla \mu (l),dl/d\tau
\right)
\end{equation}
for all $l\in \hat\mathcal{G}_{\mbox{\scriptsize red}}^{*}\subset \hat\mathcal{G}^{*}.$
Substituting now~(3.14) into~(2.8), we obtain
\begin{equation}
\{\tilde{\gamma},\tilde{\mu}\}_{{}\lambda }:=\{\gamma ,\mu \}_{{}\lambda
}\mid _{l\in \hat\mathcal{G}_{\mbox{\scriptsize red}}^{*}}=\left( \nabla
\tilde{\gamma},(\eta +\lambda \theta )\nabla \tilde{\mu}\right) ,
\end{equation}
where we have used the obvious relationship
\begin{equation}
\{\tilde{\gamma},\tilde{\mu}\}_{{}\lambda }\circ \nu =\{\tilde{\gamma}\circ
\nu ,\tilde{\mu}\circ \nu \}_{{}\lambda },
\end{equation}
and where
\be
\ba{l}
\eta :=\left(
\begin{array}{c}
 2\partial \vspace{1mm}\\
-\partial \exp [-\partial ^{-1}(\alpha +\beta )]\partial ^2-2\partial
-\partial \cdot \alpha \exp [-\partial ^{-1}(\alpha +\beta )]\cdot \partial
\end{array}
\right.  \vspace{2mm}\\
\phantom{\eta :=}\left.
\begin{array}{c}
-\partial \exp [-\partial ^{-1}(\alpha +\beta )]\partial ^2-2\partial
-\partial \cdot \beta \exp [-\partial ^{-1}(\alpha +\beta )]\cdot \partial
\vspace{1mm}\\
2\partial
\end{array}\right),
\vspace{4mm}\\
\theta :=\left(
\begin{array}{cc}
 0 & 2\partial \exp [-\partial ^{-1}(\alpha +\beta )]\partial \vspace{1mm}\\
-2\partial \exp [-\partial ^{-1}(\alpha \beta )]\partial &  0
\end{array}
\right) .
\ea
\ee
It is straightforward to verify that these integro-dif\/ferential, implectic
(=co-symplectic= Pois\-son) operators are compatible~[4] (see also~[11] for a
general theory of iso-symplectic structures on functional manifolds) on the
reduced submanifold $\hat\mathcal{G}_{\mbox{\scriptsize red}}^{*}$ and
def\/ine a bi-Hamiltonian structure on it.

\setcounter{equation}{0}

\section{Dif\/f(\textit{S}$^1$) action, associated momentum mapping \\
and versal deformations}

Let us introduce some additional notation concerning versal
deformations [1, 7]. By a de\-for\-ma\-tion of the operator~(1.1)
we shall mean an operator of the same form with a matrix
$l_{{}\lambda }(\epsilon )$ whose entries are analytic in
$\epsilon $ in a neighborhood of $\epsilon =0$ in ${\mathbb C}^n$
and satisf\/ies $l_{{}\lambda }(0)=l_{{}\lambda }$ for all
$\lambda \in {\mathbb C}.$ The coordinates $\epsilon _i\in
{\mathbb C}$, $1\leq i\leq n,$ of $\epsilon $ are called the
deformation parameters and the space of these parameters is called
the base of the deformation.

Two deformations $l_{{}\lambda }^{\prime }(\epsilon )$ and $l_{{}\lambda
}^{\prime \prime }(\epsilon )$ of a matrix $l_{{}\lambda }$ will be called
equivalent if there exists a deformation $A^{*}(\phi _{{}\epsilon
}):l_{{}\lambda }^{\prime }(\epsilon )\rightarrow l_{{}\lambda }^{\prime
\prime }(\epsilon )$ generated by a dif\/feomorphism $\phi _{{}\epsilon }\in
\, \mbox{Dif\/f}(S^1)$ satisfying $\phi _{{}\epsilon }\mid _{{}\epsilon =0}=id.$

From a given deformation $l_{{}\lambda }(\epsilon )$ one can obtain a new
deformation $\tilde{l}_{{}\lambda }(\tilde{\epsilon})$ by setting $\tilde{l}
_{{}\lambda }(\tilde{\epsilon}):=l_{{}\lambda }(\epsilon (\tilde{\epsilon})),$
where $\epsilon :{\mathbb C}^m\rightarrow {\mathbb C}^n$ is an analytic
mapping in a neighborhood of $\tilde{\epsilon}=0$ in ${\mathbb C}^m$ and
satisf\/ies the condition $\epsilon (0)=0.$ The deformation $\tilde{l}
_{{}\lambda }(\tilde{\epsilon})$ is said to be induced from $l_{{}\lambda
}(\epsilon )$ by the mapping $\epsilon :{\mathbb C}^m\rightarrow {\mathbb C}^n.$

A deformation $l_{{}\lambda }(\epsilon ),$ $\epsilon \in {\mathbb C}^n,$ is
called versal if every one of its deformations $l_{{}\lambda }(\tilde{\epsilon}),$
$\tilde{\epsilon}\in {\mathbb C}^m,$ is equivalent to a deformation
induced from it. A versal deformation is said to be universal if the induced
deformation described in the def\/inition of versality is unique.

Before we give a def\/inition of a transversal deformation for the induced
group $\hat{\mathcal{G}}_{\mbox{\scriptsize red}}$ orbits, let us consider a family of
smooth induced transformations $\phi {}_{{}\sigma }(x)\in \hat{G}_{\mbox{\scriptsize red}},$
$\sigma \in {\mathbb R},$ where $\phi _{\sigma }(x)=1+O(\sigma )$ as $\sigma
\rightarrow 0.$ Each such transformation generates (via formula (1.5)) a new
matrix $l_{{}\lambda }(\sigma ),$ $\sigma \rightarrow 0,$ that obviously
belongs to the orbit space associated to the $\hat{G}_{\mbox{\scriptsize red}}$ action. The set
of matrices
\begin{equation}
\left. \frac{dl_{{}\lambda }(\sigma )}{d\sigma }\right| _{{}\sigma =0}\in
\hat\mathcal{G}_{\mbox{\scriptsize red}}^{*}
\end{equation}
spans a linear subspace $\hat{V}_{{}\lambda }\subset
\hat\mathcal{G}_{\mbox{\scriptsize red}}^{*}$
of f\/inite codimension. Consider an arbitrary
deformation $l_{\lambda }(\epsilon ),$ $\epsilon \in {\mathbb C}^n,$ of a
given matrix $l_{\lambda }\in \hat\mathcal{G}_{\mbox{\scriptsize red}}^{*}$ and denote by
$\hat{E}_{{}\lambda }$ the linear span in
$\hat\mathcal{G}_{\mbox{\scriptsize red}}^{*}$ over the matrices $\partial
l_{{}\lambda }(\epsilon )/\partial \epsilon _i\mid _{{}\epsilon =0},$ $1\leq
i\leq n.$ The above deformation is said to be transverse to the induced
$\hat{G}_{\mbox{\scriptsize red}}$ orbit if the subspaces $\hat{E}_{{}\lambda }$ and
$\hat{V}_{{}\lambda }$ together span their ambient space, that is
\begin{equation}
\hat{E}_{{}\lambda }+\hat{V}_{{}\lambda }=\hat\mathcal{G}_{\mbox{\scriptsize red}}^{*}.
\end{equation}
The following general theorem [1] holds for versal deformations of the Dirac
operator~(1.1).

\medskip

\noindent
{\bf Theorem 4.1.} {\it A deformation $l_{{}\lambda }(\epsilon ),$
$\epsilon \in {\mathbb C}^n,$ is versal if and only if it is transverse to the
induced group $\hat{G}$ orbit.}

\medskip

This theorem can be proved by applying standard perturbation
theory techniques to the Dirac type operator~(1.1).

We are now ready to make use of the results of Section 3 to describe the
spaces $\hat{E}_{{}\lambda }$ and $\hat{V}_{{}\lambda }$ analytically. Let
$\tilde{\gamma}\in D(\hat\mathcal{G}_{\mbox{\scriptsize red}}^{*})$ be any
smooth functional on $\hat\mathcal{G}_{\mbox{\scriptsize red}}^{*};$ it
generates a f\/low on the loop group $\hat{G}_{\mbox{\scriptsize red}}$ orbit via
the $(\sigma ,x)$-evolutions
\begin{equation}
\frac{dl}{d\sigma }:=\{\tilde{\gamma},l\}_{{}\lambda },\qquad \frac{dl}{dx}:=\left(
\begin{array}{cc}
\alpha & 0 \\
0 & \beta
\end{array}
\right) l
\end{equation}
with respect to the Poisson bracket (3.16). In view of~(4.3),~(3.16) implies
that the subspace $\hat{V}_{{}\lambda }$ is isomorphic to the following
subspace of vector functions in $T^{*}(M):$
\begin{equation}
V_{{}\lambda }:=\left\{ \Lambda _{{}\lambda }\psi :=(\eta +\lambda \theta
)\psi :\nabla \tilde{\gamma}=\psi \in T^{*}(M)\right\} .
\end{equation}
Theorem 4.1 suggests the following construction of versal deformations for
the Dirac type operator~(1.1): As $\Lambda _{{}\lambda }$ is skew-symmetric,
the operator $i\Lambda _{{}\lambda }$ is formally selfadjoint in the space
$L_2({\mathbb R}/2\pi {\mathbb Z};{\mathbb C}^2).$ Therefore, the orthogonal complement
to the subspace $V_{{}\lambda }$ with respect to the natural scalar product
in $L_2({\mathbb R}/2\pi {\mathbb Z};{\mathbb C}^2)$ consists of $2\pi$-periodic
solutions to the equation
\begin{equation}
\Lambda _{{}\lambda }\psi =0.
\end{equation}
Whence we have the following characterization of versal deformations of the
operator~(1.1).

\medskip

\noindent {\bf Theorem 4.2.} {\it The prolongation of the matrix
$l_{{}\lambda }\in \hat\mathcal{G}_{\mbox{\scriptsize red}}^{*}$ defined as
\begin{equation}
\bar{l}_{{}\lambda }(\epsilon ):=\left(
\begin{array}{cc}
\lambda & \exp (\partial ^{-1}\beta ) \vspace{1mm}\\
\exp (\partial ^{-1}\alpha ) &  -\lambda
\end{array}
\right) +\sum\limits_{i,j=1}^2\epsilon _{ij}\bar{f}_i\otimes \bar{f}_j
\end{equation}
generates a versal deformation of the Dirac type operator (1.1).
Here $\otimes $ is the
usual Kronecker tensor product in ${\mathbb C}^2,$
$\epsilon_{ij}\in {\mathbb C},1\leq i,j\leq 2,$
$\epsilon _{12}=-\epsilon _{21}$ are any
deformation constants, and
$\bar{f}_i\in W_2^{(1)}({\mathbb R}/2\pi {\mathbb Z};{\mathbb C}^2),
i=1,2,$ are two linearly independent, normalized solutions to the Dirac
equations
\begin{equation}
\frac{d\bar{f}_i}{dx}+\bar{l}_{{}\lambda }\bar{f}_i=0,\qquad \left. \left\|
\bar{f}_i,\bar{f}_j\right\| \right| _{x=0}=1,
\end{equation}
with spectral parameter $\lambda \in {\mathbb C}.$}

\medskip

\noindent
{\bf Proof.} It is easy to verify that the set of solutions to
equation (4.5) is isomorphic to the set of functions
\[
\hat{\psi}=\sum\limits_{i,j=1}^2\epsilon _{ij}\bar{f}_i\otimes \bar{f}_j,
\]
and these functions satisfy the canonical Casimir equation
\begin{equation}
\left[ l_{{}\lambda },\hat{\psi}\right] -\frac{d\hat{\psi}}{dx}=0,
\end{equation}
which is equivalent to equation (4.5). Owing to the fact that any matrix
$l_{{}\lambda }\in \hat\mathcal{G}^{*}$ in~(1.1) can be
transformed into the expression $\bar{l}_{{}\lambda }(0)\in \hat\mathcal{G}^{*}$
with functional parameters $\alpha$, $\beta $
given by~(3.12), this leads to the general form~(4.6) for versal
deformations of~(1.1). This ends the proof.

\subsection*{Acknowledgment}

A. Prykarpatsky is greatful to the Dept. of Applied Mathematics at AGH
for its support of this work through an AGH research grant.
D. Blackmore would like to express his
gratitude to the Courant Institute of Mathematical Sciences for the
hospitality extended to him as a visiting member during the time when this
research was conducted.

\label{pryk-black_lp}
\end{document}